\input amstex

\mag=\magstep1\input amstex

\mag=\magstep1
\documentstyle{amsppt}
\parindent=20pt

\topmatter
\title A Base Point Free Theorem for Log canonical Surfaces
\endtitle
\author Shigetaka FUKUDA
\leftheadtext{Shigetaka FUKUDA}
\rightheadtext{Log canonical Surfaces}
\endauthor

\endtopmatter

\document

\subsubhead\nofrills{\bf Introduction}
\endsubsubhead\par
\vskip .5pc

Let $(X,\Delta)$ be a complete, log canonical algebraic surface
defined over the field of complex numbers {\bf C}.
A nef and big Cartier divisor $D$ on $X$ is {\it nef and log big} on
$(X,\Delta)$ by definition if deg$(D\vert_C)>0$ for all irreducible
components $C$ of the reduced part
$\lfloor\Delta\rfloor$ of $\Delta$.

We follow the notation and terminology of [5].

In [6] Miles Reid introduced the notion of "log big" and gave the statement as follows:

{\it Let $(M,\Gamma)$ be a complete, log canonical algebraic variety over {\bf C} and
$L$ a nef Cartier divisor on $M$. Suppose that $aL-(K_M+\Gamma)$
is nef and log big on
$(M,\Gamma)$ for some $a \in$ {\bf N}.
Then the linear system $\vert mL \vert$ is free from base points
for every $m\gg0$.}

In this paper we give a proof to this statement in the surface case.
\vskip .5pc

{\bf Main Theorem.}\quad
{\it Let $H$ be a nef Cartier divisor on $X$ such that \linebreak $aH-(K_X+\Delta)$
is nef and log big on
$(X,\Delta)$ for some $a \in$ {\bf N}.
Then the \linebreak complete linear system $\vert mH \vert$ is free from base points
for every sufficiently large integer $m$.}
\vskip .5pc

\subsubhead\nofrills{\rm REMARK 1.}
\endsubsubhead\quad
In the case where $(X,\Delta)$ is a weakly kawamata log terminal
projective surface, we gave a proof to the theorem above in [1].
\vskip .5pc

\subsubhead\nofrills{\rm REMARK 2.}
\endsubsubhead\quad
Under the assumption that $aH-(K_X+\Delta)$ is not nef and log big on
$(X,\Delta)$ but nef and big, there exists a counterexample due to
Zariski in which the theorem is not valid ([3], remark 3-1-2).
\vskip .5pc

The author would like to thank the referee for pointing out the misprints and errors in the original version of the paper.
\vskip 1pc

\subsubhead\nofrills{\bf 0. Preliminaries}
\endsubsubhead\par
\vskip .5pc

First we collect some well known results concerning normal surfaces,
which will be required for the proof of Main Theorem.

Let $\mu:V \to W$ be a birational morphism between complete,
{\bf Q}-factorial, normal algebraic surfaces over {\bf C}.
\vskip .5pc

{\bf Lemma 1 (Projection Formula).}\quad
{\it For {\bf Q}-divisors $D$ on $V$ and $G$ on $W$, \linebreak $(D,\mu^*G)=(\mu_*D,G)$.}
\vskip .5pc

{\bf Lemma 2.}\quad
{\it If $D$ is a nef {\bf Q}-divisor on $V$, $\mu_*D$ is also nef on $W$.}
\vskip .5pc

\subsubhead\nofrills{\rm Proof.}
\endsubsubhead\quad
For all irreducible curves $C$ on $W$, $(\mu_*D,C)=(D,\mu^*C) \ge 0$
from Lemma 1.
\vskip .5pc

{\bf Lemma 3.}\quad
{\it If $D$ is a big {\bf Q}-divisor on $V$, $\mu_*D$ is also big on $W$.}
\vskip .5pc

\subsubhead\nofrills{\rm Proof.}
\endsubsubhead\quad
For a Cartier divisor $A$ on $V$, $H^0(V,\Cal O_V(A))\hookrightarrow
H^0(W,\Cal O_W(\mu_{*}A))$, because $V$ and $W$ are normal.
\vskip .5pc

{\bf Lemma 4.}\quad
{\it Let $A$ be a non $\mu$-exceptional prime divisor and $B$ a nef
{\bf Q}-divisor on $V$.
$(\mu_*A,\mu_*B) \ge (A,B)$.}
\vskip .5pc

\subsubhead\nofrills{\rm Proof.}
\endsubsubhead\quad
{}From Lemma 1, $(\mu_*A,\mu_*B)=(\mu^*\mu_*A,B)$. Here $\mu^*\mu_*A
\ge A$, because $A$ is not $\mu$-exceptional. Thus $(\mu^*\mu_*A,B)
\ge (A,B)$.
\vskip .5pc

{\bf Lemma 5.}\quad
{\it Every complete, {\bf Q}-factorial, normal algebraic surface over
{\bf C} is projective.}
\vskip .5pc

\subsubhead\nofrills{\rm Proof.}
\endsubsubhead\quad
Assume that $\mu$ is a resolution of singularities of $W$ and $A$ an
ample divisor on $V$.
Then $\mu_*A$ is an ample {\bf Q}-divisor on $W$ from Lemma 3 and 4
and the Nakai-Moishezon criterion.
\vskip .5pc

Next we note a well known result concerning surface singularities, which will be used without mentioning it throughout this paper.
For the convenience of the reader we indicate a proof, which relies on the log minimal model program.
\vskip .5pc

{\bf Proposition 0.}\quad
{\it If $(X,\Delta)$ is weakly kawamata log terminal, then $X$ is {\bf Q}-factorial.}
\vskip .5pc

\subsubhead\nofrills{\rm Proof.}
\endsubsubhead\quad
Let $f: M \to X$ be a log resolution of $(X,\Delta)$ such that $K_M+f_*^{-1}\Delta+F=f^*(K_X+\Delta)+E$ with $E \geq 0$ and Supp($E$)$=$Exc($f$), where $F=\Sigma\{F_i;\quad F_i$ is an $f$-exceptional prime divisor $\}$.

Apply the relative log minimal model program to $f: (M,f_*^{-1}\Delta+F) \to X$.

We end up with a {\bf Q}-factorial weakly kawamata log terminal surface (over $X$) $g: (Y,g_*^{-1}\Delta+(F)_Y) \to X$ such that $K_Y+g_*^{-1}\Delta+(F)_Y$ is $g$-nef.

Because $(E)_Y$ is a $g$-exceptional $g$-nef divisor, $(E)_Y=0$.
Thus Exc($g$)$=\emptyset$.
Therefore $g$ is an isomorphism from Zariski's Main Theorem.
\vskip .5pc

Lastly we mention variations of results by Kawamata and Keel-Matsuki-\linebreak McKernan.
\vskip .5pc

{\bf Proposition 1.}\quad
{\it Suppose that $(X,\Delta)$ is a weakly kawamata log terminal \linebreak projective surface. Let $R$ be a $(K_X+\Delta)$-extremal
ray.
Then there exists a rational curve $C \in R$ such that
$(-(K_X+\Delta),C) \le 4$.}
\vskip .5pc

\subsubhead\nofrills{\rm Proof.}
\endsubsubhead\quad
For some $r \ge 1$, $r(K_X+\Delta)$ is Cartier.
Let $R$ be a $(K_X+\Delta)$-extremal ray.
Put $s:=$ min$\{(\Delta,E);\quad E$ is an irreducible component
of $\Delta\}$.
For $0<\varepsilon \ll \frac{1}{(\left|s\right|+1)r}$,
$K_X+(1-\varepsilon)\Delta$ is kawamata log terminal and $R$ is a
$(K_X+(1-\varepsilon)\Delta)$-extremal ray.

Thus, from [2], there exists a rational curve $C \in R$ such that
$-(K_X+(1-\varepsilon)\Delta).C \le 4$.

If $\Delta$ does not include $C$, then $(\Delta,C) \ge 0$.
Hence $-(K_X+\Delta).C \le -(K_X+(1-\varepsilon)\Delta).C \le 4$.

If $\Delta$ includes $C$, then $s \le (\Delta,C)$.
Hence
$-(K_X+\Delta).C=-(K_X+(1-\varepsilon)\Delta).C-\varepsilon(\Delta,C)
\le 4+(-s)\varepsilon$.
By the choice of $\varepsilon$, $-(K_X+\Delta).C \le 4$.
\vskip .5pc

{\bf Proposition 2.}\quad
{\it Suppose that $(X,\Delta)$ is a weakly kawamata log terminal \linebreak projective surface.
$D$ is a nef {\bf Q}-divisor, but $K_X+\Delta$ is not nef.
Set $\lambda:=$ {\rm sup}$\{\lambda \in${\bf Q}$;\quad
D+\lambda(K_X+\Delta) \quad is \quad nef\}$.

Then $\lambda$ is a rational number and moreover there is a
$(K_X+\Delta)$-extremal ray $R$ such that
$(D+\lambda(K_X+\Delta)).R=0$.}
\vskip .5pc

\subsubhead\nofrills{\rm Proof.}
\endsubsubhead\quad
{}From Proposition 1 and [4], the proof of 2.1, the assertion follows.
\vskip 1pc

\subsubhead\nofrills{\bf 1. Proof of the main theorem}
\endsubsubhead\par
\vskip .5pc

{}From a result by Shokurov ([5], 17.10) (cf. [7], 9.1) and Lemma 5, we
may find a weakly kawamata log terminal projective
surface $(Y,S+B)$ and a birational morphism $g: Y \to X$ such that
$K_Y+S+B=g^*(K_X+\Delta)$, where $S$ is the
reduced part of $S+B$.
We note that $g^*(aH-(K_X+\Delta))$ is nef and log big on
$(Y,g_*^{-1}\lfloor\Delta\rfloor+B)$.
In the case where $S=g_*^{-1}\lfloor\Delta\rfloor$, [1] implies the
assertion.
Thus we may assume that $S-g_*^{-1}\lfloor\Delta\rfloor > 0$.

We consider the following exact sequence for $m \in${\bf N}:
$$
0\to \Cal O_Y(mg^*H-(S-g_*^{-1}\lfloor\Delta\rfloor))
\to \Cal O_Y(mg^*H)
\to \Cal O_{S-g_{*}^{-1} \lfloor\Delta\rfloor}(mg^*H) \to 0
$$
Here
$mg^*H-(S-g_*^{-1}\lfloor\Delta\rfloor)-(K_Y+g_*^{-1}\lfloor\Delta
\rfloor+B) = mg^*H-(K_Y+S+B) = g^*(mH-(K_X+\Delta)) =
g^*(m-a)H+g^*(aH-(K_X+\Delta))$ is nef and log big on
$(Y,g_*^{-1}\lfloor\Delta\rfloor+B)$ for $m \ge a$.
Thus from [1],
$$
H^1(Y,\Cal O_Y(mg^*H-(S-g_*^{-1}\lfloor\Delta\rfloor)))=0
$$
Therefore the homomorphism
$$
H^0(Y,\Cal O_Y(mg^*H)) \to H^0(S-g_{*}^{-1}\lfloor\Delta\rfloor,
\Cal O_{S-g_{*}^{-1}\lfloor\Delta\rfloor}(mg^*H))
$$
is
surjective. Because dim $g(S-g_*^{-1}\lfloor\Delta\rfloor)=0$,
$$
Bs\vert mg^*H \vert \cap (S-g_*^{-1}\lfloor\Delta\rfloor) = \emptyset \tag ***
$$

Now run $(K_Y+g_*^{-1}\lfloor\Delta\rfloor+B)$-Minimal Model Program
with extremal rays that are $g^*H$-trivial (cf. [3], lemma 3-2-5 and [4]).

We have three cases:

\vskip .5pc

{\rm Case (A)}.\quad
{\it We obtain the morphism $p:Y \to Z$ such that
$g^*H=p^*(p_*g^*H)$, $p_*g^*H$ is Cartier, $(Z,p_*(g_*^{-1}\lfloor\Delta\rfloor+B))$ is a
weakly kawamata log terminal \linebreak projective surface and
$K_Z+p_*(g_*^{-1}\lfloor\Delta\rfloor+B)$ gives a non-negative
function on $\{C\in \overline{NE}(Z);(p_*g^*H,C)=0\}$.}

We put $\lambda:=$ sup$\{\lambda \in${\bf Q}$;\quad
p_*g^*H+\lambda(K_Z+p_*(g_*^{-1}\lfloor\Delta\rfloor+B))$ is nef $\}$.

If $K_Z+p_*(g_*^{-1}\lfloor\Delta\rfloor+B)$ is nef, then $\lambda=\infty$.
If $K_Z+p_*(g_*^{-1}\lfloor\Delta\rfloor+B)$ is not nef and
$\lambda=0$, then there exists a
$(K_Z+p_*(g_*^{-1}\lfloor\Delta\rfloor+B))$-extremal ray $R$ such
that $(p_*g^*H,R)=0$ from Proposition 2, but this is a contradiction
!  Thus $\lambda>0$.

We note that
$m(p_*g^*H)-p_*(S-g_*^{-1}\lfloor\Delta\rfloor)=p_*(g^*(mH)-(S-g_*^{-1}
\lfloor\Delta\rfloor))=K_Z+p_*(g_*^{-1}\lfloor\Delta\rfloor+B)+p_*(g^*(
mH)-(K_Y+S+B))=K_Z+p_*(g_*^{-1}\lfloor\Delta\rfloor+B)+p_*(g^*(mH-(K_X+
\Delta)))=K_Z+p_*(g_*^{-1}\lfloor\Delta\rfloor+B)+p_*(g^*(aH-(K_X+
\Delta)))
+p_*(g^*(m-a)H)$ is nef for $m$ such that $m-a> \frac{1}{\lambda}$.
Here $p_*(g^*(aH-(K_X+\Delta)))$ is nef and log big on
$(Z,p_*(g_*^{-1}\lfloor\Delta\rfloor+B))$ from Lemmas 2,3 and 4.

Thus $\vert m(r(p_*g^*H)-t(p_*(S-g_*^{-1}\lfloor\Delta\rfloor)))
\vert$ is base point free for $m \gg 0$ from [1], where $t$ is a
positive natural number such that
$t(p_*(S-g_*^{-1}\lfloor\Delta\rfloor))$ is
a Cartier divisor and $r$ is a sufficiently large prime number.

Therefore $Bs \vert m(p_*g^*H) \vert \subseteq
p_*(S-g_*^{-1}\lfloor\Delta\rfloor)$ for every sufficiently large
integer $m$ (cf. [3], the proof of theorem 3-1-1).
Noting the fact that $g^*H=p^*(p_*g^*H)$, we come to the conclusion
that $Bs \vert m(p_*g^*H) \vert=\emptyset$ from (***).

\vskip .5pc

{\rm Case (B)}.\quad
{\it We obtain the morphism $p:Y \to Z$, where $Z$ is a smooth curve
and $g^*H \sim p^*P$ for some divisor $P$ on $Z$.}

If deg$(P)>0$, then $\vert mP \vert $ is base point free for $m \gg 0$.
Thus $\vert mH \vert $ is base point free.

If deg$(P)=0$, then $g^*H$ is numerically trivial.
{}From (***), $\vert g^*(mH) \vert \ne \emptyset $ for $m \ge a$.
Thus $mH$ is linearly trivial.

\vskip .5pc

{\rm Case (C)}.\quad
{\it We obtain the morphism $p:Y \to Z$, where $Z$ is a point
 and $g^*H$ is linearly trivial.}
\vskip .5pc

\centerline{$\hbox to 7cm{\hrulefill}$}

\par\bigskip\parindent=100pt

Shigetaka Fukuda

Suzuka National College of Technology (Suzuka K{\B o}sen)

JAPAN

\bigskip

Current address:

Faculty of Education, Gifu Shotoku Gakuen University

Yanaizu-ch{\B o}, Gifu City, Gifu Prefecture

501-6194 JAPAN

e-mail:\quad fukuda\@ha.shotoku.ac.jp
  
\end{document}